\documentclass[11pt, leqno]{article}
\usepackage{amsmath, amsfonts}
\begin{document}
\author{Ajai Choudhry}
\title{Symmetric diophantine systems and\\ families of  elliptic curves of high rank}
\date{}
\maketitle

\begin{abstract} 
While there has been considerable interest in the problem of finding elliptic curves of high rank over $\mathbb{Q}$, very few parametrized families of elliptic curves of generic rank $\geq 8$ have been published. In this paper we use solutions of certain symmetric diophantine systems to construct several parametrized families of elliptic curves with their generic ranks ranging from at least 8 to at least 12. Specific numerical values of the parameters yield elliptic curves with quite large coefficients and we could therefore determine the precise rank only in a few cases where the rank of the elliptic curve $\leq 13$. It is, however, expected that the parametrized families of elliptic curves obtained in this paper would yield examples of elliptic curves with much higher rank. 
\end{abstract}



\section{Introduction}

This paper is concerned with elliptic curves defined over the field $\mathbb{Q}$ of rational numbers. According to Mordell's celebrated theorem on elliptic curves, the set of rational points $E(\mathbb{Q})$ on an elliptic curve defined over $\mathbb{Q}$ may be given the structure of a finitely generated  abelian group, and we may write
$E(\mathbb{Q}) \simeq T \oplus Z^r$  where $T$ is the torsion subgroup and $r$ is a nonnegative integer known as the rank of the elliptic curve.

The vast majority of elliptic curves have rank  $ \leq 2$. In fact, of all the elliptic curves with conductor $< 10^8$ given in the Stein-Watkins database of elliptic curves, $98.1\%$  have rank $\leq 2$ \cite[p. 251]{BMSW}. Further, it has been stated by Silverman \cite[p. 254]{Si} that ``it is difficult to produce curves $E/\mathbb{Q}$ having even moderately high rank".   There has thus been considerable interest in constructing elliptic curves of large rank. At present, the largest known rank of an individual elliptic curve is 28 \cite{El}. It has been conjectured that there is no upper bound on $r$ \cite[p. 254]{Si} but there is some  evidence to suggest that probably there are only finitely many elliptic curves with rank exceeding 21 (\cite{PP}, \cite{Po}). 

In this paper we construct a number of  families of elliptic curves whose coefficients are given in terms of several arbitrary rational parameters and whose generic rank ranges from at least 8 to at least 12. In this context, it is pertinent to note that while several authors have constructed  parametrized families of elliptic curves of rank $\leq 6$, very few parametrized families of elliptic curves with generic rank $\geq 8$ have been published. Two families of elliptic curves of generic rank 8 in terms of one parameter and five parameters have been given  by Shioda \cite{Sh} and by Fermigier \cite{Fe} respectively. One-parameter families of elliptic curves of generic ranks $\geq 11$ and $\geq 12$ were constructed by Mestre (\cite{Me1}, \cite{Me2}), and using Mestre's method, Nagao \cite{Na} found a one-parameter family of rank $\geq 13$. Subsequently, Kihara \cite{Ki} gave a one-parameter family of elliptic curves of rank $\geq 14$.

We will use solutions of certain symmetric diophantine systems to construct our examples of families of parametrized curves of high rank. The manner of construction ensures that we will know a certain number of rational points on these elliptic curves. We illustrate  the general approach  of producing  these families in Section 2 constructing inter alia two families of elliptic curves of rank 4 and 7. In Section 3 we construct parametrized  families of elliptic curves of higher ranks. 

For specific numerical values of the parameters, our parametrized families  yield examples of elliptic curves over $\mathbb{Q}$ with quite large  coefficients. This is not surprising since, as has been pointed out by Cassels \cite[p. 257]{Ca1}, ``the theory makes it clear that an abelian variety can only have high rank if it is defined by equations with very large coefficients". The large numerical coefficients were a major obstacle both in conducting efficient searches for specific elliptic curves of high rank and in determining the ranks of individual elliptic curves. We could accordingly  determine the precise rank only in a few cases where the rank of the elliptic curve $\leq 13$.

\section{A description  of the  method of constructing \\parametrized families of elliptic curves}

Let us assume that there exists an identity,
\begin{equation}
(x+a_1)(x+a_2)(x+a_3)+a_0^2=(x+b_1)(x+b_2)(x+b_3)+b_0^2, \label{identcub1}
\end{equation}
where $x$ is arbitrary and $a_i,\,b_i,\;i=0,\,1,\,2,\,3$, are rational numbers. If we denote either 
side of the identity \eqref{identcub1} by $\phi(x)$, it immediately follows from \eqref{identcub1} that $\phi(x)$ becomes a perfect square when $x$ is assigned any of the 6 values $a_i,\,b_i,\;i=1,\,2,\,3$, and hence we know  6 rational points on the cubic curve $y^2=\phi(x)$.  

To construct the identity \eqref{identcub1}, we need to solve the simultaneous diophantine equations,
\begin{equation}
\begin{aligned}
a_1+a_2+a_3&= b_1+b_2+b_3,\\
a_1a_2+a_2a_3+a_3a_1&=b_1b_2+b_2b_3+b_3b_1,\\
a_1a_2a_3+a_0^2&=b_1b_2b_3+b_0^2.
\end{aligned} \label{diophsys1}
\end{equation}
In view of the relations between the elementary symmetric functions and sums of powers, the first two equations of the diophantine system \eqref{diophsys1} are equivalent to the Tarry-Escott problem of degree 2 whose complete solution is well-known \cite[p. 52]{Di}. Using this solution, the third equation is readily solved since any rational number can be expressed as the difference of the  squares of two rational numbers.

Since the aforementioned complete solution of the Tarry-Escott problem of degree 2 is expressed in  terms of 5 arbitrary parameters, we have  obtained a cubic curve $y^2=\phi(x)$ whose coefficients are  in terms of 5 parameters and on which we already know 6 rational points. In general, the discriminant of $\phi(x)$ is not 0, and hence we have constructed a parametrized family of  elliptic curves on which we know  6 rational points. With reference to the group of rational points $E(\mathbb{Q})$, since the three  points $(-a_i,\,a_0),\;i=1,\,2,\,3$ on our elliptic curve lie on the straight line $y=a_0$, at most  two of these points can be independent points of infinite order. Similarly, at most  two of the points $(-b_i,\,b_0)$ can be independent. We should  expect  the generic rank of our parametrized family of elliptic curves  to be 4.  Indeed, the solution, 
\[(a_0,\,a_1,\,a_2,\,a_3,\,b_0,\,b_1,\,b_2,\,b_3)=(96,\,10,\,13,\,  -23,\,66,\,-17,\, -5,\,  22),\]  
of the diophantine system \eqref{diophsys1} yields the elliptic curve, 
\begin{equation}
y^2=x^3-399x+6226, \label{ecex1}
\end{equation}
 with 4 of the known points, namely $(-10,\;96),\;(-13,\;96),\;(5,\;66),\;(17,\;66)$, being independent points of infinite order.   The rank of the elliptic curve \eqref{ecex1} is, in fact,  4. It follows that the generic rank of our parametrized family of elliptic curves is at least 4.

Similarly, instead of the identity \eqref{identcub1}, we may construct an identity,
 \begin{equation}
(x+a_1)(x+a_2)(x+a_3)(x+a_4)+a_0^2=(x+b_1)(x+b_2)(x+b_3)(x+b_4)+b_0^2, \label{identquart1}
\end{equation}
by using the complete solution of the Tarry-Escott problem \cite[p. 55-58]{Di} and thus obtain a parametrized family of quartic curves 
$y^2=\phi(x)$  on which we know  8 rational points, their  abscissae  being $-a_i,\,-b_i,\;i=1,\,\ldots,\,4$. In general, the quartic curve $y^2=\phi(x)$ represents a quartic model of an elliptic curve and following a well-known procedure (see, for instance, Mordell \cite[p. 77]{Mo}, or Cassels \cite[pp. 35-36]{Ca2}), we can find a birational transformation that would reduce the quartic model $y^2=\phi(x)$ to the cubic Weierstrass model of an elliptic curve. We thus obtain a parametrized family of elliptic curves on which we know 8  rational points. 

As an example, when we take 
\begin{equation*}
\begin{aligned}
(a_0,\,a_1,\,a_2,\,a_3,\,a_4)&=(9333234,\, -1940, \, -1076, \, 1324, \, 1692),\\
(b_0,\,b_1,\,b_2,\,b_3,\,b_4)&=(9541134,\, -2196,\, -356,\, 460,\, 2092),
\end{aligned}
\end{equation*}
the relation \eqref{identquart1} is an identity, and we get the quartic model of an elliptic curve,
\begin{equation}
y^2=x^4-4768608x^2-460748288x+91785556686276, \label{ecex2}
\end{equation}
with the 8 known rational points on it being $(-a_i,\,a_0),\;(-b_i,\,b_0),\;i=1,\ldots,\,4$. Eq.~\eqref{ecex2} may be reduced by a birational transformation to the Weierstrass form of the elliptic curve given by
\begin{equation}
 y^2 = x^3 - 23420131301937 + 18114867816009096080, \label{ecex2red}
\end{equation}
and corresponding to the 8 known points on the quartic curve \eqref{ecex2}, we get 8 rational points on the elliptic curve \eqref{ecex2red}. It has been verified that 7 of these 8 points are independent points of infinite order, and hence it follows that the generic rank of the parametrized family of elliptic curves that we obtain from the identity \eqref{identquart1} is at least 7.
 
More generally, in the identities \eqref{identcub1} and \eqref{identquart1}, $a_0$ and $b_0$ could be taken as suitable polynomials of degree $\leq  2$ in the variable $x$, and we can still obtain elliptic curves on which we know  6 or 8 rational points. With such minor modifications in the procedure outlined above, we can construct several parametrized families of elliptic curves whose generic rank will, in general,  be at most 8. 

We can generalise the above procedure even further and  construct identities of the type,
\begin{equation}
f_1^2(x)-\Pi_{i=1}^{2n}(x+a_i)=f_2^2(x)-\Pi_{i=1}^{2n}(x+b_i), \label{identgen}
\end{equation}
where $a_i,\,b_i,\;i=1,\,2,\,\ldots,\,2n$, are rational numbers and $f_1(x),\;f_2(x)$ are suitably chosen polynomials of degree $n$ such that each side of the identity \eqref{identgen} reduces to a cubic or a quartic polynomial. To construct such identities, we have to solve certain symmetric diophantine systems. We will always denote the cubic or quartic polynomial on either side of an identity  of the type \eqref{identgen}  by $\phi(x)$. We thus obtain parametrized families of elliptic curves $y^2=\phi(x)$ on which we know $4n$ rational points, their abscissae being $-a_i,\;-b_i,\;i=1,\ldots,\,2n$.  
 
Once  we have constructed a  family of elliptic curves with coefficients given in terms of several parameters, say $t_1,\,t_2,\,\ldots,\,t_r$ and we know certain rational points on these elliptic curves, we need to check whether these are  independent points of infinite order. In view of a well-known theorem of Silverman, 
to prove the independence of a set of  points $P_j,\;j=1,\,2,\,\ldots,\,m$, expressed in parametric terms, it suffices to show that for a certain specialization $t_1^{\prime},\,t_2^{\prime},\,\ldots,\,t_r^{\prime}$ of rational values of the parameters $t_1,\,t_2,\,\ldots,\,t_r$, the points $P_j,\;j=1,\,2,\,\ldots,\,m$, are independent over the corresponding elliptic curve defined over $\mathbb{Q}$. For this purpose, we invoke another well-known  theorem \cite[Theorem 8.1, p. 242]{SZ} according to which,  if the regulator of the points $P_j,\;j=1,\,2,\,\ldots,\,m$, is not zero, the points must be independent points of infinite order. 

Having constructed a parametrized family of elliptic curves, it is natural to try to find specific rational values of the parameters  which yield examples of  elliptic curves of high rank. As the numerical size of the coefficients in most of our examples is just too large, it was, in general,  not feasible to follow the usual search strategies (eg by computing  Mestre-Nagao sums)  and we had  to necessarily  restrict ourselves to small values of the parameters to obtain our numerical examples.

\section{Parametrized families of elliptic curves}
In this section  we use solutions of certain symmetric diophantine systems to  obtain several  parametrized families  of elliptic curves with generic ranks ranging from  at least 8 to at least 12. 

\subsection{}
We will  construct our first example of a parametrized family of elliptic curves by using solutions of the following simultaneous symmetric diophantine chains:
\begin{equation}
\begin{aligned}
a_1+a_2+a_3&=&b_1+b_2+b_3&=&c_1+c_2+c_3&=&d_1+d_2+d_3&=0,\\
\psi(a_1,\,a_2,\,a_3)&=&\psi(b_1,\,b_2,\,b_3)&=&\psi(c_1,\,c_2,\,c_3)&=&\psi(d_1,\,d_2,\,d_3)&=s,
\end{aligned}
\label{tepdeg2chns}
\end{equation}
where $\psi(t_1,\,t_2,\,t_3)=t_1t_2+t_2t_3+t_3t_1$. 

If we write
\begin{equation}
\begin{aligned}
f_1(u_1,\,u_2,\,v_1,\,v_2,\,w_1,\,w_2)&=u_1v_1w_1-u_1v_2w_2-u_2v_1w_2\\
& \quad \quad \quad -u_2v_2w_1-u_2v_2w_2,\\
f_2(u_1,\,u_2,\,v_1,\,v_2,\,w_1,\,w_2)&=-u_1v_1w_1-u_1v_1w_2-u_1v_2w_1\\
& \quad \quad \quad -u_2v_1w_1+u_2v_2w_2,
\end{aligned}
\label{deffuvw}
\end{equation}
a readily verifiable solution of the simultaneous diophantine chains \eqref{tepdeg2chns} is given by
\begin{equation}
\begin{aligned}
a_1&=f_1(p_1,\,p_2,\,q_1,\,q_2,\,r_1,\,r_2),\;&a_2&=f_2(p_1,\,p_2,\,q_1,\,q_2,\,r_1,\,r_2),\\
b_1&=f_1(p_2,\,p_1,\,q_1,\,q_2,\,r_1,\,r_2),\;&b_2&=f_2(p_2,\,p_1,\,q_1,\,q_2,\,r_1,\,r_2),\\
c_1&=f_1(p_1,\,p_2,\,q_2,\,q_1,\,r_1,\,r_2),\;&c_2&=f_2(p_1,\,p_2,\,q_2,\,q_1,\,r_1,\,r_2),\\
d_1&=f_1(p_1,\,p_2,\,q_1,\,q_2,\,r_2,\,r_1),\;&d_2&=f_2(p_1,\,p_2,\,q_1,\,q_2,\,r_2,\,r_1),\\
a_3&=-a_1-a_2,\;\;b_3=-b_1-b_2,&c_3&=-c_1-c_2,\;\;d_3=-d_1-d_2,
\end{aligned}
\label{soltepdeg2chns}
\end{equation}
and
\begin{equation}
s=-(p_1^2+p_1p_2+p_2^2)(q_1^2+q_1q_2+q_2^2)(r_1^2+r_1r_2+r_2^2),
\end{equation}
where $p_1,\,p_2,\,q_1,\,q_2,\,r_1,\,r_2$ are arbitrary  parameters.

 For the sake of brevity, we will write,
\begin{equation}
\pi_1=a_1a_2a_3,\;\;\pi_2=b_1b_2b_3,\;\;\pi_3=c_1c_2c_3,\;\;\pi_4=d_1d_2d_3.
\end{equation}

When $a_i,\,b_i,\,c_i,\,d_i,\;i=1,\,2,\,3$, are  given by \eqref{soltepdeg2chns}, the  relations \eqref{tepdeg2chns} are satisfied and we  have,
\begin{multline}
\Pi_{i=1}^3(x+a_i)\Pi_{i=1}^3(x+b_i)-\Pi_{i=1}^3(x+c_i)\Pi_{i=1}^3(x+d_i)\\
=(x^3+sx+\pi_1)(x^3+sx+\pi_2)-(x^3+sx+\pi_3)(x^3+sx+\pi_4)\\
=(\pi_1+\pi_2-\pi_3-\pi_4)(x^3+sx)+\pi_1\pi_2-\pi_3\pi_4\quad \quad \quad \quad \quad  \\
=\{x^3+sx+(\pi_1\pi_2-\pi_3\pi_4)/(\pi_1+\pi_2-\pi_3-\pi_4)+(\pi_1+\pi_2-\pi_3-\pi_4)/4\}^2\\
\quad -\{x^3+sx+(\pi_1\pi_2-\pi_3\pi_4)/(\pi_1+\pi_2-\pi_3-\pi_4)-(\pi_1+\pi_2-\pi_3-\pi_4)/4\}^2,
\end{multline}
and, on transposition, we get the identity,
\begin{multline}
\{x^3+sx+(\pi_1\pi_2-\pi_3\pi_4)/(\pi_1+\pi_2-\pi_3-\pi_4)\\
+(\pi_1+\pi_2-\pi_3-\pi_4)/4\}^2-\Pi_{i=1}^3(x+a_i)\Pi_{i=1}^3(x+b_i)\\
=\{x^3+sx+(\pi_1\pi_2-\pi_3\pi_4)/(\pi_1+\pi_2-\pi_3-\pi_4)\\
-(\pi_1+\pi_2-\pi_3-\pi_4)/4\}^2-\Pi_{i=1}^3(x+c_i)\Pi_{i=1}^3(x+d_i). \label{ident1}
\end{multline}
Each side of the identity \eqref{ident1} is, in fact, a cubic polynomial in the variable $x$. If we denote this cubic polynomial by $\phi(x)$,  it is clear from \eqref{ident1} that $\phi(x)$  becomes a perfect square when $x$ takes any of the 12 rational values $-a_i,\,-b_i,\,-c_i,\,-d_i,\;i=1,2,3$. Further, the discriminant of 
$\phi(x)$ is, in general, nonzero, and hence the cubic equation $y^2=\phi(x)$ represents an elliptic curve on which we know 12 rational points,  their abscissae being  $-a_i,\,-b_i,\,-c_i,\,-d_i,\;i=1,2,3$.

On making the invertible linear transformation given by
\begin{equation}
\begin{aligned}
x&=-X/\{2(\pi_1+\pi_2-\pi_3-\pi_4)(\pi_1-\pi_2-\pi_3+\pi_4)(\pi_1-\pi_2+\pi_3-\pi_4)\},\\
 y& = -Y/\{4(\pi_1+\pi_2-\pi_3-\pi_4)^2(\pi_1-\pi_2-\pi_3+\pi_4)(\pi_1-\pi_2+\pi_3-\pi_4)\},
\end{aligned}
\end{equation}
the elliptic curve $y^2=\phi(x)$ may be written, after suitable simplification, as follows: 
\begin{multline}
Y^2= X^3+(\pi_1-\pi_2-\pi_3+\pi_4)^2(\pi_1-\pi_2+\pi_3-\pi_4)^2(\pi_1+\pi_2-\pi_3-\pi_4)^2\\
\times \left[ 4sX+\left \{2\sum_{i=1}^4\pi_i^2-\left (\sum_{i=1}^4\pi_i \right)^2\right \}^2-64\pi_1\pi_2\pi_3\pi_4\right ]\quad \quad \quad \quad \label{ecfamilytepdeg2XY}
\end{multline}
Further, we know  12 rational points on the elliptic curve \eqref{ecfamilytepdeg2XY}, their    abscissae being  given by $ka_i,\, kb_i,\,kc_i,\,kd_i,\;i=1,\,2,\,3$ where 
\[k=2(\pi_1+\pi_2-\pi_3-\pi_4)(\pi_1-\pi_2-\pi_3+\pi_4)(\pi_1-\pi_2+\pi_3-\pi_4).\]

It follows from \eqref{ident1} that $\phi(-a_1)=\phi(-a_2)=\phi(-a_3)$ and hence, of the three points on the curve $y^2=\phi(x)$ with abscissae $-a_i,\;i=1,\,2,\,3$, only two can be independent. A similar remark is applicable to the triples of points with abscissae $\{-b_i\},\;\{-c_i\},\;\{-d_i\}$ and hence  at most  8 of the  12 known rational points on the elliptic curve \eqref{ecfamilytepdeg2XY} can be independent.

As a numerical example, when we take $p_1=5,\, p_2= -1,\, q_1=   7,\, q_2=   1,\, r_1=   11,\, r_2=   15$, the elliptic curve \eqref{ecfamilytepdeg2XY} may be written, after suitable reduction, as follows:
\begin{equation}
Y^2 = X^3-22064074044012X+43046837966291058900.
\label{ecfamilytepdeg2XYex1a}
\end{equation}

The following 8 points out of the 12 known rational points on the curve \eqref{ecfamilytepdeg2XYex1a} are independent points of infinite order:
\begin{equation}
\begin{aligned}
(2648646, 1785505722),\quad  & (2774772, 1785505722),\\
 (3441438, 2805949146),\quad &  (1909908, 2805949146), \\
	(-1639638, -8649618978),\quad  & (5297292, -8649618978), \\
  (-702702, -7629175554),\quad  &(5009004, -7629175554).
	\end{aligned} \label{ecfam2ex1points}
	\end{equation}

Using the software SAGE,  we  found the following three additional points  on the curve \eqref{ecfamilytepdeg2XYex1a}: 
\begin{equation}
\begin{aligned}
(-667436, 7581284326), \quad & (-426972, 7238075754), \\
(430353, 5799241371). \quad &
\end{aligned}
\label{ecfam2ex1points2}
\end{equation}

The regulator of the  11  points given by \eqref{ecfam2ex1points} and \eqref{ecfam2ex1points2}, as  computed by the software SAGE, is 1381592532.65. Thus these 11 points are independent and the rank of the elliptic curve \eqref{ecfamilytepdeg2XYex1a} is at least 11. In fact, according to APECS (a package written in MAPLE for working with elliptic curves), the rank of the curve \eqref{ecfamilytepdeg2XYex1a} is $\leq 11$. Since we already have 11 independent points on the curve \eqref{ecfamilytepdeg2XYex1a}, its rank must be 11. Further, the torsion group of the curve \eqref{ecfamilytepdeg2XYex1a} is trivial. 

Since in the special case of  the  elliptic curve \eqref{ecfamilytepdeg2XYex1a}, we found 8 independent points of infinite order out of the 12 known rational points,   it follows that  the generic rank of the family of elliptic curves given  by  \eqref{ecfamilytepdeg2XY} is at least 8.

\subsection{}

The family of elliptic curves given in this subsection was obtained following the general approach described in Section 2 but is presented below in a simpler, albeit somewhat different, way.

	Let there exist  rational numbers $a_i,\,b_i,\,c_i,\;i=1,\,2,\,3,\,4$, such that 
\begin{equation}
\sum_{i=1}^4a_i^r=\sum_{i=1}^4b_i^r=\sum_{i=1}^4c_i^r,\quad r=1,\,2,\,3. \label{tepdeg3chn1}
\end{equation}
We then have the identities,
\begin{align}
\Pi_{i=1}^4(x-a_i)-\Pi_{i=1}^4(x-b_i)&=a_1a_2a_3a_4-b_1b_2b_3b_4,\label{ident1ab}\\
\Pi_{i=1}^4(x-a_i)-\Pi_{i=1}^4(x-c_i)&=a_1a_2a_3a_4-c_1c_2c_3c_4.\label{ident1ac}
\end{align}

We now consider the quartic curve defined by
\begin{multline}
y^2=8(2a_1a_2a_3a_4-b_1b_2b_3b_4-c_1c_2c_3c_4)(x-a_1)(x-a_2)\\
\times (x-a_3)(x-a_4)+(b_1b_2b_3b_4-c_1c_2c_3c_4)^2.\quad \quad \quad \quad \quad \quad \label{ecfamily1}
\end{multline}

It is clear that four rational points on the  curve \eqref{ecfamily1} are given by $(a_i,\,b_1b_2b_3b_4-c_1c_2c_3c_4), \;i=1,\,2,\,3,\,4$. We also note that in view of \eqref{ident1ab}, we have
\begin{equation*}
\begin{aligned}
y^2&=8(2a_1a_2a_3a_4-b_1b_2b_3b_4-c_1c_2c_3c_4)\\
& \quad \quad \times \{(x-b_1)(x-b_2)(x-b_3)(x-b_4)+a_1a_2a_3a_4-b_1b_2b_3b_4\}\\
& \quad \quad \quad +(b_1b_2b_3b_4-c_1c_2c_3c_4)^2,\\
&=8(2a_1a_2a_3a_4-b_1b_2b_3b_4-c_1c_2c_3c_4)(x-b_1)(x-b_2)\\
& \quad \quad \times (x-b_3)(x-b_4)+(4a_1a_2a_3a_4-3b_1b_2b_3b_4-c_1c_2c_3c_4)^2.
\end{aligned}
\end{equation*}
It follows that the four rational points $(b_i, \, 4a_1a_2a_3a_4-3b_1b_2b_3b_4-c_1c_2c_3c_4)$,\; $i=1,\,2,\,3,\,4$ also lie on the elliptic curve \eqref{ecfamily1}. 

Similarly, in view of the identity \eqref{ident1ac}, we may write \eqref{ecfamily1} equivalently as 
\begin{multline*}
y^2=8(2a_1a_2a_3a_4-b_1b_2b_3b_4-c_1c_2c_3c_4)(x-c_1)(x-c_2)\\
\times (x-c_3)(x-c_4) +(4a_1a_2a_3a_4-b_1b_2b_3b_4-3c_1c_2c_3c_4)^2.
\end{multline*}
and hence we have four more points $(c_i,\,4a_1a_2a_3a_4-b_1b_2b_3b_4-3c_1c_2c_3c_4)$,\;  $i=1,\,2,\,3,\,4$ on the curve \eqref{ecfamily1}.

We note that, in general, Eq.~\eqref{ecfamily1} represents a quartic model of an elliptic curve, and we already know  12 rational points on this  curve. We will use one of the 12 known rational points to  find a birational transformation that would reduce the quartic model \eqref{ecfamily1} to the cubic Weierstrass model of an elliptic curve, and corresponding to the 11 remaining known rational points on the quartic curve \eqref{ecfamily1}, we will get  11 rational points on the cubic model of the elliptic curve.  

We will now obtain solutions of the multigrade chains \eqref{tepdeg3chn1}. We note that  solutions of \eqref{tepdeg3chn1} obtained by writing $a_3=-a_1,\,a_4=-a_2,\,b_3=-b_1,$ $ b_4=-b_2,\,c_3=-c_1,\,c_4=-c_2$, and thereby reducing the system of equations \eqref{tepdeg3chn1} to the diophantine chain $a_1^2+a_2^2=b_1^2+b_2^2=c_1^2+c_2^2$ lead to examples of elliptic curves on which at most 6 of the above 12 points are independent. We therefore exclude such solutions of \eqref{tepdeg3chn1} from consideration.

 It has been shown in \cite{ChW} that if there exists a solution  of the diophantine chains,
\begin{equation}
\begin{aligned}
x_1^2+y_1^2+z_1^2&=x_2^2+y_2^2+z_2^2&=x_3^2+y_3^2+z_3^2,\\
x_1y_1z_1&=x_2y_2z_2&=x_3y_3z_3,\quad \quad \;\;
\end{aligned}
\label{eqsumsprodchn}
\end{equation}
and we write
\begin{equation}
\begin{aligned}
a_1& = x_1 - y_1 - z_1,\quad \;\;& a_2& = y_1 - z_1 - x_1,\;\;\\
 a_3& = z_1 - x_1 - y_1,\;\; & a_4& = x_1 + y_1 + z_1,\\
b_1& = x_2 - y_2 - z_2,& b_2& = y_2 - z_2 - x_2, \\
b_3& = z_2 - x_2 - y_2, & b_4& = x_2 + y_2 + z_2,\\
c_1& = x_3 - y_3 - z_3,& c_2& = y_3 - z_3 - x_3, \\
c_3& = z_3 - x_3 - y_3, & c_4& = x_3 + y_3 + z_3,
\end{aligned} \label{valabc}
\end{equation}
then the numbers $a_i,\, b_i, \,c_i,\;i=1,\,2,\,3,\,4$, satisfy the diophantine chains \eqref{tepdeg3chn1}. Several  parametric solutions of the simultaneous diophantine equations \eqref{eqsumsprodchn} are given in \cite{ChW}. These solutions yield parametrized families of the cubic model of elliptic curves  with 11 known rational points. 

As an example, the following one-parameter solution (see \cite[p.\ 93]{ChW}) of the simultaneous diophantine chains \eqref{eqsumsprodchn}, 
\begin{equation} 
\begin{aligned}
x_1 &= 2(p^2+1)^3,\\
y_1 &= -(p - 1)(p + 1)(p^2-2p-1)(p^2+2p-1),\\
z_1 &=  2p(p^2-2p-1)(p^2+2p-1),\\
x_2 &= (p^2+1)(p^2+2p-1)^2,\\
y_2 &= -2(p-1)(p^2+1)(p^2-2p-1),\\
z_2 &= 2p(p+1)(p^2+1)(p^2-2p-1),\\ 
x_3 &=(p^2+1)(p^2-2p-1)^2,\\
y_3 &=  2(p+1)(p^2+1)(p^2+2p-1),\\
z_3 &= -2p(p-1)(p^2+1)(p^2+2p-1),
\end{aligned}
\end{equation}
where $p$ is an arbitrary parameter, yields the following parametrized family of elliptic curves: 
\begin{multline}
y^2= x^3-27(11p^{16}+56p^{14}-12p^{12}+1672p^{10}-638p^8+1672p^6-12p^4+56p^2+11)\\
\times (251p^{40}+216p^{38}+32658p^{36}+216984p^{34}+671007p^{32}+2551136p^{30}\\
-5621544p^{28}+42013152p^{26}-67017498p^{24}+222343248p^{22}-127186644p^{20}\\
+222343248p^{18}-67017498p^{16}+42013152p^{14}-5621544p^{12}+2551136p^{10}\\
+671007p^8+216984p^6+32658p^4+216p^2+251)x+54(p^2+1)^2\\
\times (11p^{16}+56p^{14}-12p^{12}+1672p^{10}-638p^8+1672p^6-12p^4+56p^2+11)^2\\
\times (1195p^{48}-2398p^{46}+170128p^{44}+1528142p^{42}+8633422p^{40}+17765158p^{38}\\
+68866832p^{36}+222561418p^{34}-1080958651p^{32}+5138210548p^{30}-7489078176p^{28}\\
+18959483244p^{26}-11645588604p^{24}+18959483244p^{22}-7489078176p^{20}\\
+5138210548p^{18}-1080958651p^{16}+222561418p^{14}+68866832p^{12}+17765158p^{10}\\
+8633422p^8+1528142p^6+170128p^4-2398p^2+1195). \label{cubicecfamily1}
\end{multline} 

There are 11 known rational points on this elliptic curve whose co-ordinates, in terms of the parameter $p$,  are too cumbersome to write and therefore,  we do not give these 11 points explicitly. It has been verified that for several values of $p$, 9 out of the 11 known rational points are independent. 

As a numerical example, when $p=2$, we get the elliptic curve,
\begin{multline}
y^2=x^3-1136261157571019728659411x\\
+466685506089477132791551368731316450,\quad \quad \quad \quad \quad \quad \label{ecfamily1ex1} 
\end{multline}
on which, out of the set of 11 known rational points,  there are 9 independent points. 
The invertible linear transformation,
\begin{equation}
x = 1089X+363,\quad  y = 35937Y,
\end{equation}
reduces the elliptic curve \eqref{ecfamily1ex1} to the minimal Weierstrass form given below:
\begin{multline}
Y^2=X^3+X^2-958125505468762024X\\
+361360495869188941993242116. \quad \quad \quad \quad \quad \quad  \label{ecfamily1ex1W}
\end{multline}

The 9 independent points on the curve \eqref{ecfamily1ex1W}  corresponding to  the  9 independents on the curve \eqref{ecfamily1ex1} are as follows:
\begin{equation}
\begin{aligned}
&(11149617229984/19321, 2123284512273843750/2685619), \\
&(-1171372325542/2401, 3140764291484625000/117649), \\
&(751945051961/16, 651907077554613387/64), \\
&(49012446586/81, 1294374710916872/729), \\
&(198014368166/361, 6305852501956812/6859), \\
&(641648531, 3279574736568), \\
&(59134648706/121, 4174284741150888/1331), \\
&(3005025814/9, 239485247913716/27), \\
&(9096019456, 862686952661118).
\end{aligned}
\label{ecfamily1ex1Wpoints1}
\end{equation}

We could find the following 4 additional  independent points on the curve \eqref{ecfamily1ex1W} using SAGE:
\begin{equation}
\begin{aligned} (-48407290, \,20189782661724),\;\;& (14065754, 18651716303160),\\
(73997912, 17054805737250),\;\;& (90932456, 16582752629118).\
\end{aligned} \label{ecfamily1ex1Wpoints2}
\end{equation}

The regulator of the 13 points given by \eqref{ecfamily1ex1Wpoints1} and \eqref{ecfamily1ex1Wpoints2}, as computed by SAGE, is 389828159565.83 confirming that these 13 points are independent. According to APECS, the upper limit for the rank of the curve \eqref{ecfamily1ex1W} is 13, and since we have already found 13 independent points on this curve, its rank must be 13.

Since   9 out of the 11 known  points on the curve \eqref{cubicecfamily1} are independent in the special case when $p=2$,  the generic rank of the family of elliptic curves given  by  \eqref{cubicecfamily1} is at least 9. 

As a second example,  when $p=5$, the minimal form of the elliptic curve \eqref{cubicecfamily1} is given by
\begin{multline}
y^2=x^3+x^2-7526665968750696273075230520630360x\\
+249763118169814320753109971139754041307175559002084. \label{cubicecfamily1ex2}
\end{multline}
The rank of the curve \eqref{cubicecfamily1ex2} was determined by the software MAGMA to be 13. The abscissae of the 13 independent points of infinite order on the curve \eqref{cubicecfamily1ex2} are given below:
\begin{equation*}
\begin{aligned}
&48634078038299278376983377/76247824,\\&
8070419987630413400822239075641585/61155659809263616,\\&
6092992423962972607322092830628647/80977443542798881,\\&
461089280807619481351249719657621/8249557371401476,\\&
9289851754473511254669830882814460671/35684558324214797449,\\&
52525257722608835220923418493449617799/49899550759585419649,\\&
277522049066244995131051766066574/4107959693396689,\\&
20640417956268123829481926974/159239161,\\&
2109895471283944241447243629999998/32199609729595009,\\&
2510782654688490433512918809080605/16428907017388804,\\&
86250313881356948400248550959543664531/1211480391364699110769,\\&
1749307288076044116212179533960194190/9270856201405973329,\\&
3892604721565311997214458198413742127151/  \\
& \quad 16881711960419335668409.
\end{aligned}
\end{equation*}

When $p=7$, the curve \eqref{cubicecfamily1} is given by
\begin{multline}
y^2=x^3-70080465281896617533991165744653594072570451x+\\
226274285487270710244792369861166227599729743261105332834795161250. \label{cubicecfamily1ex3}
\end{multline}
It has been verified using SAGE that 9 out of the 11 known points on this curve are independent points and hence the rank of the curve \eqref{cubicecfamily1ex3} is  $\geq 9$. The upper bound for the analytic rank of the curve \eqref{cubicecfamily1ex3} as determined by SAGE is 20. We could not determine the precise  rank because of the size of the coefficients. 

It has been shown in \cite{ChW} that infinitely many three-parameter solutions of the diophantine chains \eqref{eqsumsprodchn} can be obtained. Each of these solutions yields a  three-parameter solution of the diophantine chains \eqref{tepdeg3chn1}, and hence we can obtain infinitely many  families of elliptic curves whose coefficients are given in terms of three arbitrary parameters. As the three-parameter solutions of \eqref{eqsumsprodchn} are very cumbersome to write, these families of elliptic curves are even more cumbersome. We note that elliptic curves belonging to these parametrized families  are expected to have high rank. 

Finally we note that by interchanging the quadruples $\{a_i\},\;\{b_i\},\;\{c_i\}$, we can get two more parametrized families of elliptic curves whose generic rank will be at least 9.

\subsection{}
In this section we will construct a parametrized family of elliptic curves using two sets of rational numbers  $a_i,\;b_i,\;i=1,\,2,\,\ldots,\,6$, satisfying certain conditions. We will denote the elementary symmetric functions of the six rational numbers $a_i$ by $s_i,\;i=1,\,2,\,\ldots,\,6$, that is,
\begin{equation*}
s_1=\sum_{i=1}^6a_i,\quad s_2=\sum_{i<j}a_ia_j, \quad s_3=\sum_{i<j<k}a_ia_ja_k,\quad s_6=a_1a_2\cdots a_6,
\end{equation*}
and similarly we will denote the corresponding elementary symmetric functions of the rational numbers $b_i$ by $t_i,\;i=1,\,2,\,\ldots,\,6$. We  also write $d_i=s_i-t_i,\;i=1,\,2,\,\ldots,\,6$. We will use this notation  in the next two subsections also.

Let  there exist rational numbers  $a_i,\;i=1,\,2,\,\ldots,\,6$, and $b_i,\;i=1,\,2,\,\ldots,$ $ 6$, satisfying the symmetric diophantine system,
\begin{equation}
s_1=t_1, \quad s_2=t_2, \quad s_1s_3-2s_4=t_1t_3-2t_4. \label{diophsysfam3}
\end{equation}

We then have
\begin{multline}
\Pi_{i=1}^6(x+a_i)-\Pi_{i=1}^6(x+b_i)=d_3x^3+d_4x^2+d_5x+d_6\\
=\{x^3+(d_4x^2+d_5x+d_6)/d_3+d_3/4\}^2-\{x^3+(d_4x^2+d_5x+d_6)/d_3-d_3/4\}^2,\label{ident1fam3}
\end{multline}
and we get the identity,
\begin{multline}
\{x^3+(d_4x^2+d_5x+d_6)/d_3+d_3/4\}^2-\Pi_{i=1}^6(x+a_i)\\
\quad =\{x^3+(d_4x^2+d_5x+d_6)/d_3-d_3/4\}^2-\Pi_{i=1}^6(x+b_i).
\label{ident1afam3}
\end{multline}

Denoting either side of the identity \eqref{ident1afam3} by $\phi(x)$, and using the relations \eqref{diophsysfam3}, we get,
\begin{multline}
\phi(x)=-(d_3^2s_2-2d_3d_5-d_4^2)x^4/d_3^2\\
+(d_3^3-2d_3^2s_3+4d_3d_6+4d_4d_5)x^3/(2d_3^2)\\
+(d_3^2d_4-2d_3^2s_4+4d_4d_6+2d_5^2)x^2/(2d_3^2)+(d_3^2d_5-2d_3^2s_5+4d_5d_6)x/(2d_3^2)\\
+(d_3^4+8d_3^2d_6-16d_3^2s_6+16d_6^2)/(16d_3^2) \quad \quad \quad \label{defphifam3}
\end{multline}

Now $\phi(x)$ is, in general, a quartic polynomial in $x$ and it follows from \eqref{ident1afam3} that $\phi(x)$ becomes a perfect square when $x$ is assigned any of the 12 values $-a_i, \, -b_i, \;i=1,\,2,\,\ldots,\,6$. We therefore know  12 rational points, with abscissae $-a_i,\,-b_i,\;i=1,\,2,\,\ldots,\,6$, on the quartic curve
\begin{equation}
y^2=\phi(x), \label{ecfam3}
\end{equation}
where $\phi(x)$ is defined by \eqref{defphifam3}.

We will now obtain solutions of the symmetric diophantine system \eqref{diophsysfam3}. Since our objective is to generate a family of elliptic curves, we will exclude those  solutions of Eqs.~\eqref{diophsysfam3} for which  the discriminant of $\phi(x)$ becomes 0. 

A simple solution of the diophantine system \eqref{diophsysfam3} is given by 
\begin{equation}
\begin{aligned}
a_1& = (p+2q)m-(p-q)n, &a_2 &= -(2p+q)m-(p+2q)n,\\
 a_3& = (p-q)m+(2p+q)n, &a_4& = (r+2s)u-(r-s)v, \\
a_5 &= -(2r+s)u-(r+2s)v, &a_6& = (r-s)u+(2r+s)v,\\
b_1 &= (p-q)m-(p+2q)n, &b_2 &= -(2p+q)m-(p-q)n,\\
b_3 &= (p+2q)m+(2p+q)n, &b_4 &= (r-s)u-(r+2s)v,\\
 b_5 &= -(2r+s)u-(r-s)v, &b_6 &= (r+2s)u+(2r+s)v,
\end{aligned} \label{ecfam3sol1}
\end{equation}
where $m,\,n,\,p,\,q,\,r,\,s,\,u,\,v$ are arbitrary parameters.

When the values of $a_i,\,b_i$ are given by \eqref{ecfam3sol1}, the discriminant of  $\phi(x)$ is, in general, not zero, and hence Eq.~\eqref{ecfam3} represents a quartic model of an elliptic curve on which we know  12 rational points whose abscissae are given by $-a_i, \, -b_i,\;i=1,\,2,\,\ldots,\,6$.  Using any of these known points, we can reduce the quartic equation \eqref{ecfam3} to a cubic model with 11 known rational points. 

As a numerical example, when $(m,\,n,\,p,\,q,\,r,\,s,\,u,\,v)=(1,\, 7,\, 3,\, 5,\, 2,\,  11,$ $ 3,\,4)$, 
we get the elliptic curve,
\begin{equation}
Y^2 + Y = X^3 - 78654091314536101X + 4993138309311379361023650. \label{ecfam3ex1}
\end{equation}
Nine of the 11 known rational points on the curve \eqref{ecfam3ex1} are independent. These points are given below:
\begin{equation}
\begin{aligned}
&(1419528264, 52476440000169),\\
&(15427375761/16, -1838910202387815/64),\\
&(12233281459/25, -1143413169906927/125),\\
&(144154185/16, 132485022151835/64),\\
&(52499795215/841, 13933569409840090/24389),\\
&(-153755798, -3667661902713),\\
&(-2478831521/25, -429700761309558/125),\\
&(1680856689/25, 11737369812738/125),\\
&(-512313311/9, -82276801860745/27).
\end{aligned}
\label{pointsecfam3ex1}
\end{equation}

Using SAGE, we found the following additional independent point on the curve \eqref{ecfam3ex1}:
\begin{equation}
(29506697560865/68644, -127942115051260580605/17984728).\label{pointaddlfam3ex1}
\end{equation}
The regulator of the 10 points \eqref{pointsecfam3ex1} and \eqref{pointaddlfam3ex1}, as computed by SAGE, is 44242748.70, confirming that these points on the curve \eqref{ecfam3ex1} are independent. In fact, it appears that the rank of the elliptic curve \eqref{ecfam3ex1} is 10.

 Since 9 out of the 11 known rational points on the curve \eqref{ecfam3ex1} are independent, it follows that the generic rank of the parametrized family of elliptic curves given by \eqref{ecfam3} is at least 9. 

We will now obtain a parametrized family of elliptic curves of generic rank at least 10 by obtaining another parametric solution of the symmetric diophantine system \eqref{diophsysfam3}. 

To solve the diophantine system \eqref{diophsysfam3}, we write,
\begin{equation}
\begin{aligned}
a_1 &= mp_1v+nv, \\
a_2 &= mp_2v-\{(p_1-p_2)u+v\}n, \\
a_3 &= mp_3v+\{(p_1-2p_2+p_3)u+v\}n,\\
 a_4 &= mp_4v-\{(p_1-2p_2+2p_3-p_4)u+v\}n, \\
a_5 &= mp_5v+\{(p_1-2p_2+2p_3-2p_4+p_5)u+v\}n,\\
 a_6 &= (p_1-p_2+p_3-p_4+p_5)mv+\{(p_2-p_3+p_4-p_5)u-v\}n,\\
 b_1 &= mp_1v-nv, \\
b_2 &= mp_2v+\{(p_1-p_2)u+v\}n, \\
b_3 &= mp_3v-\{(p_1-2p_2+p_3)u+v\}n, \\
b_4 &= mp_4v+\{(p_1-2p_2+2p_3-p_4)u+v\}n,\\
 b_5 &= mp_5v-\{(p_1-2p_2+2p_3-2p_4+p_5)u+v\}n,\\
 b_6 &= (p_1-p_2+p_3-p_4+p_5)mv-\{(p_2-p_3+p_4-p_5)u-v\}n,
\end{aligned}
\label{subs2ecfam3}
\end{equation}
where $m,\,n,\,u,\,v$ and $p_i,\;i=1,\ldots,\,5$ are arbitrary parameters. With these values of $a_i, \, b_i$, 
 it is readily verified that $s_1=t_1$ and $s_2=t_2$. Further,
\begin{multline*}
s_1s_3-2s_4-(t_1t_3-2t_4)=4mnv(m^2v^2-n^2u^2)
\{f_1(p_1,\,p_2,\,p_3,\,p_4,\,p_5)u\\+f_2(p_1,\,p_2,\,p_3,\,p_4,\,p_5)v\}, \quad \quad \quad \quad 
\end{multline*}
where
\begin{equation}
\begin{aligned}
&f_1(p_1,\,p_2,\,p_3,\,p_4,\,p_5)=(p_2p_3-2p_2p_4+p_2p_5-p_3^2+p_3p_4\\
& \quad+p_4p_5-p_5^2)p_1^2-(3p_2^2p_3-6p_2^2p_4+3p_2^2p_5-4p_2p_3^2+8p_2p_3p_4\\
& \quad-4p_2p_3p_5-4p_2p_4^2+8p_2p_4p_5-4p_2p_5^2+p_3^3-2p_3^2p_4+p_3^2p_5\\
& \quad+p_3p_4^2-4p_3p_4p_5+3p_3p_5^2+3p_4^2p_5-4p_4p_5^2+p_5^3)p_1\\
& \quad+(2p_3-4p_4+2p_5)p_2^3-3(p_3-2p_4+p_5)(p_3-p_4+p_5)p_2^2\\
& \quad+(p_3^3-6p_3^2p_4+5p_3^2p_5+7p_3p_4^2-12p_3p_4p_5+5p_3p_5^2-2p_4^3+7p_4^2p_5\\
& \quad-6p_4p_5^2+p_5^3)p_2+(p_4-p_5)(p_3-p_4)(p_3^2+2p_3p_5-2p_4p_5+p_5^2),\\
&f_2(p_1,\,p_2,\,p_3,\,p_4,\,p_5)=(p_2-p_3+p_4-p_5)p_1^2-(p_2^2-2p_2p_3+4p_2p_4\\
& \quad-2p_2p_5+p_3^2-2p_3p_4+p_4^2-2p_4p_5+p_5^2)p_1-(p_3-3p_4+p_5)p_2^2\\
& \quad+(p_3-3p_4+p_5)(p_3-p_4+p_5)p_2+(p_3-p_4)(p_3+p_5)(p_4-p_5).
\end{aligned}
\end{equation}

We accordingly take
\begin{equation}
u=f_2(p_1,\,p_2,\,p_3,\,p_4,\,p_5),\quad v=-f_1(p_1,\,p_2,\,p_3,\,p_4,\,p_5), \label{valuvfam3}
\end{equation}
and a solution of the simultaneous equations \eqref{diophsysfam3}, in terms of arbitrary parameters $m,\,n,\,p_1,\,p_2,\,\ldots,\,p_5$, is given by \eqref{subs2ecfam3} where $u$ and $v$ are given by \eqref{valuvfam3}.

The above solution of the diophantine system \eqref{diophsysfam3} yields  a second family of elliptic curves \eqref{ecfam3} whose coefficients 
 are given in terms of arbitrary parameters $m,\,n,\,p_1,\,p_2,\,\ldots,\,p_5$, and on which we know 12 rational points whose abscissae are $-a_i, \, -b_i,\;i=1,\,2,\,\ldots,\,6$, where the values of $a_i,\,b_i$ are given by \eqref{subs2ecfam3}.

As a numerical example, when we take  \[(p_1,\,p_2,\,p_3,\,p_4,\,p_5,\,m,\,n)=(17,\,13,\,3,\,4,\,5,\,7,\,11),\]  
we get, on appropriate scaling,  the following solution of the diophantine system \eqref{diophsysfam3}:
\begin{equation}
\begin{aligned}
(a_1,\,a_2,\,a_3,\,a_4,\,a_5,\,a_6)&=(780,\,1228,\,1314,\, -1207,\, 1398,\,-1413),\\
 (b_1,\,b_2,\,b_3,\,b_4,\,b_5,\,b_6)&=(648, \,-136,\, -1062,\,  1543,\,  -978,\,  2085),
\end{aligned}
\label{solabex2fam3}
\end{equation}
and we thus obtain the quartic curve,
\begin{multline}
y^2=4834227853x^4+2608103498364x^3-10273867756103916x^2\\
-3211739790886292400x+6210528187439825204004, \label{quartex2fam3}
\end{multline}
on which we know 12 known rational points with abscissae $-a_i,\,-b_i$, where $a_i,\,b_i$ are given by \eqref{solabex2fam3}. Using the first known point $(-780,\,54925180902)$, we transform the curve \eqref{quartex2fam3} to the following elliptic curve:
\begin{multline}
Y^2=X^3-10228318920208466879353470719100X\\
+12588116153737599336213325703794679840709018225. \label{cubicex2fam3}
\end{multline}

We now know  11  rational points on the elliptic curve \eqref{cubicex2fam3} which are as follows:
\begin{equation*}
\begin{aligned}
&(344439380978543286825/50176, \\
& \quad 5795192172625238725997832515115/11239424), \\ 
&(37864559525453094264/7921, 190378377726108145643306415087/704969), \\
&(7390525362088551811530/3948169, \\
& \quad  7350096643832089074996920406555/7845011803), \\ 
&(40752177956740667940/10609,  189196926839822274401570687145/1092727), \\  
&(1880956160305134, 1977925012570157318577), \\
&(815024696079665265/4,  -735703752927007217019962295/8), \\  
&(348393114269455662474/52441, \\
& \quad  -5856860499598881854668281177693/12008989), \\  
&(176126802575085971520/94249,  -5850422935694673933554684595/28934443), \\
&(1815574360830758787030/582169, \\
& \quad  46632563221988934899907398562195/444194947), \\
&(1879426025359431, -1826271996173185607046), \\
&(51312686300223412434/21025,  141754394796246133687133659377/3048625).
\end{aligned}
\end{equation*}
The regulator of the first ten of these points, as computed by SAGE, is 2078733082632.16 which shows that they are independent. It follows that the generic rank of the family of the curves \eqref{ecfam3} is at least 10.

We can, in fact, choose the parameters $p_1,\,p_2,\,p_3,\,p_4,\,p_5$, such that the coefficient of $x^4$ in the quartic curve \eqref{ecfam3} is a perfect square, and then we will have 12 known rational points on the corresponding cubic model of the elliptic curve. As a numerical example, when we take 
\[(p_1,\,p_2,\,p_3,\,p_4,\,p_5,\,m,\,n)=(2,\,25/24,\,1/12,\,1,\,553/24,\,6,\,169),\]
and proceed as in the above example, we get, after appropriate scaling, the quartic curve,
\begin{multline}
y^2=26^2x^4+1678476x^3-35730288264x^2\\
+910919663233536x+7876475744393292201,
\end{multline}
which reduces to the cubic elliptic curve,
\begin{multline}
y^2 + y = x^3 - 1262161549165629265725x \\
 + 17269863219444345499893710734056, \quad \quad \quad \quad  \label{ecfam3ex3}
\end{multline}
on which we know  12 rational points. However, only ten of these points are independent. This reconfirms that the generic rank of the family of elliptic curves \eqref{ecfam3} is at least 10. We could not determine the precise rank of the elliptic curve \eqref{ecfam3ex3}.

\subsection{}
Let  there exist two sets of rational numbers $a_i,\,i=1,\,2,\,\ldots,\,6$, and $b_i,\,i=1,\,2,\,\ldots,\,6$, satisfying the symmetric diophantine system,
\begin{equation}
s_1=t_1=0,\quad s_3=t_3, \quad s_6=t_6. \label{diophsysfam4}
\end{equation}

We then have
\begin{multline}
\Pi_{i=1}^6(x+a_i)-\Pi_{i=1}^6(x+b_i)=d_2x^4+d_4x^2+d_5x\\
=\{x^3+(d_4x+d_5)/d_2+d_2x/4\}^2-\{x^3+(d_4x+d_5)/d_2-d_2x/4\}^2,\label{ident1fam4}
\end{multline}
and we get the identity,
\begin{multline}
\{x^3+(d_4x+d_5)/d_2+d_2x/4\}^2-\Pi_{i=1}^6(x+a_i)\\
=\{x^3+(d_4x+d_5)/d_2-d_2x/4\}^2-\Pi_{i=1}^6(x+b_i).\label{ident1afam4}
\end{multline}

Denoting either side of the identity \eqref{ident1afam4} by $\phi(x)$, and using the relations \eqref{diophsysfam4}, we get,
\begin{multline}
\phi(x)=\{-8((s_2+t_2)d_2-4d_4)d_2x^4-16(d_2s_3-2d_5)d_2x^3\\
+(d_2^4-8(s_4+t_4)d_2^2+16d_4^2)x^2 \quad \\
-8((s_5+t_5)d_2^2-4d_4d_5)x-16d_2^2s_6+16d_5^2\}/(16d_2^2). \label{defphifam4}
\end{multline}

Now $\phi(x)$ is a quartic polynomial in $x$ and it follows from \eqref{ident1afam4} that $\phi(x)$ becomes a perfect square when $x$ is assigned any of the 12 values $-a_i, \, -b_i, i=1,\,2,\,\ldots,\,6$.

It is readily verified that a solution of the symmetric diophantine system \eqref{diophsysfam4} is given by 
\begin{equation}
\begin{aligned}
a_1 &= h_1p(pq-r^2),  & a_2 &= -h_1q(p^2-qr),  & a_3 &= h_1r(pr-q^2),   \\
a_4 &= h_2u(uv-w^2),  & a_5 &= -h_2v(u^2-vw),  & a_6 &= h_2w(uw-v^2),  \\
b_1 &= h_1q(pq-r^2),  & b_2 &= -h_1r(p^2-qr),  & b_3 &= h_1p(pr-q^2),  \\
 b_4 &= h_2v(uv-w^2),  & b_5 &= -h_2w(u^2-vw),  & b_6 &= h_2u(uw-v^2),  
\end{aligned} \label{valabfam4}
\end{equation}
where $p,\,q,\,r,\,u,\,v,\,w,\,h_1$ and $h_2$ are arbitrary parameters. Further, with these values of $a_i, \,b_i$, the discriminant of the polynomial $\phi(x)$ is, in general, nonzero and hence the equation
\begin{equation}
y^2=\phi(x), \label{ecfam4}
\end{equation}
where $\phi(x)$ is defined by \eqref{defphifam4}, represents a quartic model of an elliptic curve and we know 12 rational points on the curve \eqref{ecfam4}, their  abscissae being $-a_i,\,-b_i,\,i=1,\,2,\,\ldots,\,6$.

When the values of $a_i, \,b_i$ are given by \eqref{valabfam4}, the coefficients of the quartic polynomial $\phi(x)$ are too cumbersome to write and hence we do not give them explicitly. We will restrict ourselves to a couple of specific examples.

As a numerical example, when we take \[(p,\,q,\,r,\,u,\,v,\,w,\,h_1,\,h_2)=(1,\, 2,\,   3 ,\,  1 ,\,  2,\,   6 ,\,  1,\,   1),\]  we get the  quartic model of an elliptic curve,
\begin{multline}
y^2=228921813x^4+825628494x^3+40295401173x^2\\
-65980539336x+183592277856,\quad \quad \quad \quad \quad \quad \label{ecfam4ex1}
\end{multline}
on which we know 12 rational points. Using the standard method already mentioned earlier, we reduce the curve \eqref{ecfam4ex1} to the  Weierstrass model,
\begin{multline}
Y^2+XY=X^3-X^2-47739269184667111896X\\
+34123010787688902778640228336, \quad \quad \quad \quad \quad \quad \label{ecfam4ex1cubic}
\end{multline}
on which we know  11  rational points. We verified that the following 8 points out of the 11 known rational points are independent:
\begin{equation}
\begin{aligned}
&(15442990559, 1726208463100830),\\
&(8073744052, -418307928289508),\\
&(56888714071, 13469519053795315),\\
&(903891853220/121, 408863422459744464/1331),\\
&(6694247596, -120558361349036),\\
&(1610576393500/169, 1465076329698280768/2197),\\
&(7007739736, 209084125480324),\\
&(29343760275580/4489, 11163685217533742272/300763).
\end{aligned} \label{pointsecfam4ex1a}
\end{equation}

Using SAGE, we found the following three additional independents on the  elliptic curve \eqref{ecfam4ex1cubic}:
\begin{equation}
\begin{aligned}
&(75830334587992/6241, -550761451920442092356/493039),\\
&(1811433219931/100, 2261228723749097149/1000),\\
&(179007778397401/5041, 2350163270404407688064/357911).
\end{aligned} \label{pointsecfam4ex1b}
\end{equation}

The regulator of the  11 points on the curve \eqref{ecfam4ex1cubic} given by \eqref{pointsecfam4ex1a} and \eqref{pointsecfam4ex1b}, as computed by SAGE, is  78091770934.92, thus confirming the independence of these points. It appears that the rank of the elliptic curve \eqref{ecfam4ex1cubic} is, in fact,  11. 

We note that if we take 
\begin{equation}
\begin{aligned} h_1&= (u-v)^2(v-w)^2(u-w)^2(u+v+w)^2(uv+uw+vw)^2\\
& \quad \times (pq-r^2)(p^2-qr)(pr-q^2)pqr,\\
 h_2 &=(p-q)^2(q-r)^2(p-r)^2(p+q+r)^2(pq+pr+qr)^2\\
 & \quad \times (u^2-vw)(uv-w^2)(uw-v^2)uvw,
\end{aligned}
\label{valhfam4}
\end{equation}
 the constant term on the right-hand side of \eqref{defphifam4} vanishes and on applying the transformation $y=Y/X^2,\,x=1/X$,  the quartic model \eqref{ecfam4} of the elliptic curve reduces to the  cubic model and corresponding to the 12 known rational points on the quartic model of the curve, we now get 12 rational points on the cubic model of the elliptic curve. 

As a numerical example, when  we  take $(p,\,q,\,r,\,u,\,v,\,w)=(2,\,3,\,7,\,1,\,$ $ 6,\,5)$,   we get the elliptic curve,
\begin{multline*}
y^2=x^3 - 1087315540936651467921677192224906456/1011981676729x^2\\
 - 2080666617378148813778824077800168117359619615625000000/1005973x\\
 + 3550434424189702535356505800072940355252343426161415124149\\
138805670907141898925975000000000/1018026243284102317, \label{ecfam4ex2}
\end{multline*}
on which we know  12 rational points. Nine of these 12 points are independent, their abscissae being 
\begin{equation*}
\begin{aligned}
&-1559069322059990914623052102500/1005973, \\ 
&2629018856807043503089852565000/1005973, \\
&3830856048490263390216642309000/1005973,\\ 
&2114036106613885301522619945000/1005973, \\
&-230843023135998969706492982500/1005973, \\
&-1039379548039993943082034735000/1005973,\\ 
&1126722367203018644181365385000/1005973, \\
&352339351102314216920436657500/1005973, \\
&-277011627763198763647791579000/1005973.
\end{aligned}
\end{equation*}
The  regulator of the above 9 points, as computed by the software SAGE, is 1769919385554.43. Thus, the generic rank of the family of elliptic curves \eqref{ecfam4} is at least 9.

\subsection{}
We will construct our final example of a family of elliptic curves using  two sets of rational numbers $a_i,\,i=1,\,2,\,\ldots,\,6$, and $b_i,\,i=1,\,2,\,\ldots,\,6$, satisfying the symmetric diophantine system,
\begin{equation}
s_1=t_1=0,\quad s_3=t_3, \quad s_2s_3-2s_5=t_2t_3-2t_5. \label{diophsysfam5}
\end{equation}

If $h$ is any arbitrary rational number, we have 
\begin{equation}
\begin{aligned}
&\{\Pi_{i=1}^6(x+a_i)-\Pi_{i=1}^6(x+b_i)\}(x^2-h^2)\\
& \quad =(d_2x^4+d_4x^2+d_5x+d_6)(x^2-h^2)\\
& \quad =\{x^4+(d_4x^2+d_5x+d_6)/d_2+d_2(x^2-h^2)/4\}^2\\
&\quad \;\;-\{x^4+(d_4x^2+d_5x+d_6)/d_2-d_2(x^2-h^2)/4\}^2, \label{ident1fam5}
\end{aligned}
\end{equation}
and we thus get the identity,
\begin{multline}
\{x^4+(d_4x^2+d_5x+d_6)/d_2+d_2(x^2-h^2)/4\}^2\\
-(x^2-h^2)\Pi_{i=1}^6(x+a_i)=\{x^4+(d_4x^2+d_5x+d_6)/d_2 \quad \quad \quad \\
-d_2(x^2-h^2)/4\}^2-(x^2-h^2)\Pi_{i=1}^6(x+b_i). \label{ident2fam5}
\end{multline}

Denoting the polynomial on either side of the identity \eqref{ident2fam5} by $\phi(x)$, we note that $\phi(x)$ is a sextic polynomial in which   the coefficients of $x^6$ and $x^5$ are given  by
\begin{equation}
(2d_2h^2+d_2^2-2d_2s_2+4d_4)/(2d_2)\;\;\; {\rm and}\;\;\;   -(d_2s_3-2d_5)/d_2
\end{equation}
respectively. In view of the relations \eqref{diophsysfam5}, the coefficient of $x^5$ is  easily seen to be 0. 

We will choose $a_i,\,b_i$ and $h$ such that the coefficient of $x^6$ also becomes 0, that is, the numbers $a_i,\,b_i$ and $h$ must satisfy the diophantine equation,
\begin{equation}
2d_2h^2+d_2^2-2d_2s_2+4d_4=0. \label{cfx6fam5} 
\end{equation}
 This  reduces $\phi(x)$ to a quartic polynomial that may be written as follows:
\begin{multline}
\phi(x)=\left[\{8(s_2+t_2)d_2^2h^2+d_2^4-8(s_4+t_4)d_2^2+32d_2d_6+16d_4^2\}x^4 \right.\\
+\{16d_2^2h^2s_3-8(s_5+t_5)d_2^2+32d_4d_5\}x^3\\
-\{2h^2d_2^4-8((s_4+t_4)h^2-s_6-t_6)d_2^2-32d_4d_6-16d_5^2\}x^2\\
+\left. \{8(s_5+t_5)d_2^2h^2+32d_5d_6\}x+d_2^4h^4+8(s_6+t_6)d_2^2h^2+16d_6^2\right]/(16d_2^2).
\end{multline}

It  readily follows  from the identity  \eqref{ident2fam5} that $\phi(-a_i),\; \phi(-b_i),\; i=1,\,\,2,$  $\ldots,\,6$, as well as $\phi(h)$ and $ \phi(-h)$, are all perfect squares. Accordingly,  we know 14 rational solutions of the quartic equation
\begin{equation}
y^2=\phi(x) \label{ecfam5}.
\end{equation}

We will now solve the symmetric diophantine system \eqref{diophsysfam5}
together with the diophantine equation \eqref{cfx6fam5}. For this purpose, 
we write
\begin{equation}
\begin{aligned}
a_1 &= (p_1+2q_1)m_1-(p_1-q_1)n_1+r, \\a_2 &= -(2p_1+q_1)m_1-(p_1+2q_1)n_1+r,\\
 a_3 &= (p_1-q_1)m_1+(2p_1+q_1)n_1+r, \\a_4 &= -(p_1-q_1)m_1+(p_1+2q_1)n_1-r,\\
 a_5 &= (2p_1+q_1)m_1+(p_1-q_1)n_1-r, \\a_6 &= -(p_1+2q_1)m_1-(2p_1+q_1)n_1-r,\\
b_1 &= (p_2+2q_2)m_2-(p_2-q_2)n_2+r, \\b_2 &= -(2p_2+q_2)m_2-(p_2+2q_2)n_2+r,\\
 b_3 &= (p_2-q_2)m_2+(2p_2+q_2)n_2+r, \\b_4 &= -(p_2-q_2)m_2+(p_2+2q_2)n_2-r,\\
 b_5 &= (2p_2+q_2)m_2+(p_2-q_2)n_2-r, \\b_6 &= -(p_2+2q_2)m_2-(2p_2+q_2)n_2-r, 
\end{aligned}
\label{valabfam5}
\end{equation}
where $p_1,\,p_2,\,q_1,\,q_2,m_1,\,m_2,\,n_1,\,n_2$ and $r$ are arbitrary nonzero parameters. It is readily seen that with the values of $a_i,\,b_i$, given by \eqref{valabfam5}, $s_1=t_1=0$ while both of the remaining two equations of \eqref{diophsysfam5} reduce to the following equation:
\begin{equation}
m_1n_1p_1q_1(p_1+q_1)(m_1+n_1) = m_2n_2p_2q_2(p_2+q_2)(m_2+n_2). \label{eqmnfam5}
\end{equation}
We now write $m_2=u_1m_1,n_2=u_2n_1$ when Eq.~\eqref{eqmnfam5} reduces to a linear equation in $m_1$ and $n_1$  and hence we get the following complete solution of Eq.~\eqref{eqmnfam5}:
\begin{equation}
\begin{aligned}
m_1 &= p_1^2q_1+p_1q_1^2-u_1u_2^2p_2^2q_2-u_1u_2^2p_2q_2^2,\\
 n_1 &= -p_1^2q_1-p_1q_1^2+u_1^2u_2p_2^2q_2+u_1^2u_2p_2q_2^2,\\
m_2 &= u_1(p_1^2q_1+p_1q_1^2-u_1u_2^2p_2^2q_2-u_1u_2^2p_2q_2^2), \\
n_2 &= -u_2(p_1^2q_1+p_1q_1^2-u_1^2u_2p_2^2q_2-u_1^2u_2p_2q_2^2),
\end{aligned} \label{valmnfam5}
\end{equation}
where $p_1,\,p_2,\,q_1,\,q_2,\,u_1,\,u_2$ are arbitrary nonzero parameters.

To solve \eqref{cfx6fam5}, we write for brevity,
\begin{equation}
\begin{aligned}
f_1(p_i,\,q_i,\,m_i,\,n_i)&=(p_1^2+p_1q_1+q_1^2)(m_1^2+m_1n_1+n_1^2)\\
& \quad \quad -(p_2^2+p_2q_2+q_2^2)(m_2^2+m_2n_2+n_2^2),\\
f_2(p_i,\,q_i,\,m_i,\,n_i)&=(p_1-q_1)(p_1+2q_1)(2p_1+q_1)(m_1-n_1)\\
& \quad \quad \times (m_1+2n_1)(2m_1+n_1) -(p_2-q_2)(p_2+2q_2)\\
& \quad \quad \times (2p_2+q_2)(m_2-n_2)(m_2+2n_2)(2m_2+n_2),
\end{aligned}
\end{equation}
and now, with the values of $a_i,\,b_i$ given by \eqref{valabfam5}, we may write Eq.~\eqref{cfx6fam5} as follows:
\begin{equation}
f_1(p_i,\,q_i,\,m_i,\,n_i)(h^2+3r^2)-f_2(p_i,\,q_i,\,m_i,\,n_i)r=0. \label{cfx6afam5}
\end{equation}
On writing, $h=vr$, we immediately get the following solution of Eq.~\eqref{cfx6afam5}:
\begin{equation}
\begin{aligned}
h&=f_2(p_i,\,q_i,\,m_i,\,n_i)v/\{(v^2+3)f_1(p_i,\,q_i,\,m_i,\,n_i)\},\\
r&=f_2(p_i,\,q_i,\,m_i,\,n_i)/\{(v^2+3)f_1(p_i,\,q_i,\,m_i,\,n_i)\},
\end{aligned}
\label{valhrfam5}
\end{equation}
where $v$ is an arbitrary parameter.

Thus, a solution of the simultaneous diophantine equations \eqref{diophsysfam5} and \eqref{cfx6fam5} is given by \eqref{valabfam5} where 
the values of $m_1,\,m_2,\,n_1,\,n_2$ are given by \eqref{valmnfam5} and the values of $h$ and $r$ are given by \eqref{valhrfam5} while   
$p_1,\,p_2,\,q_1,\,q_2,\,u_1,\,u_2$ and $v$ are arbitrary rational parameters.

With the values $a_i,\,b_i$ and $h$ chosen as above, the discriminant of the quartic polynomial $\phi(x)$ is, in general, nonzero and hence Eq.~\eqref{ecfam5} represents the  quartic model of a parametrized family of  elliptic curves on which we know $14$ rational points whose abscissae are given by  $a_i, \,b_i, i=1,\,2,\,\ldots,\,6$ and $\pm h$.

Using any of the known rational points and following the usual procedure, we can find a birational transformation that would reduce the above family of curves to  a family of elliptic curves in Weierstrass form  on which we will have 13 known rational points. As both the quartic and the Weierstrass models of this family of elliptic curves are too cumbersome to write, we do not explicitly present  the quartic and the cubic equations of this family but restrict ourselves to giving a couple of numerical examples.

We examined a  number of elliptic curves belonging to the family given by \eqref{ecfam5}. In almost all of these cases,   at least 11 of the 13 known rational points were found to be independent and   in some few cases we found 12 points to be independent. We could not find any example in which all of the 13 points were independent. 

As a numerical example,  when we take the values of the parameters $p_1,\,p_2,\,q_1,\,q_2,\,u_1,\,u_2$ and $v$ as $1,\,3,\,1,\,2,\,-1,\,-2$ and 2 respectively, we get an elliptic curve whose Weierstrass  is as follows:
\begin{multline}
y^2 = x^3-1355064646307559724826793297x\\
+19084107576037868085853647087238447797889. \quad \quad \quad \label{fam5ex1}
\end{multline}
The following 11  points on this curve are independent:
\begin{align*}
&(18418864238101437/4489, 35065928670711263766652015/300763), \\
& (3515162745891513/256, -226819310502255350099885/4096), \\
&  (16221384779701687/9, 2065579906231919277602675/27), \\
& (2864636680992817/144, 346463627191057660375/1728), \\
&  (14187172304236697/4, -1689742260465066483837525/8), \\
& (41558257019243, 185862260604690940725), \\
& (111473132761217/4, 435573324632366466105/8), \\
& (19724021739113, -5491624688580711135), \\
& (1246644171412377/64, -4573407885714980840125/512),\\
& (8631074783686853287/199809, \\
& \quad \quad 18118450186381927653890152625/89314623), \\
& (15822599097429157563/1125721, \\
& \quad \quad -63367327980237642659414589655/1194389981). \\
\end{align*}

The regulator of the   11 points given above, as computed by SAGE, is  1379591921192.48. It follows that the rank of the elliptic curve \eqref{fam5ex1} is at least 11. 
We could neither  determine any additional independent points on the curve \eqref{fam5ex1} nor its rank on account of the size of the coefficients.

As a second example, when we take the values of the parameters $p_1,\,p_2,$ $q_1,\,q_2,\,u_1,\,u_2,\,v$  as $1,\,3,\,1,\,2,\,-1,\,7,\,10$ respectively, we get the    elliptic curve, 
\begin{scriptsize}
\begin{multline} 
y^2 = x^3-1297408166125308307368483393939351175899224072475192243x\\+537211144047153364926448357254872490627314988850400537154755515461210216146647438,
\end{multline}
\end{scriptsize}
on which we have the following 12 points:
\begin{scriptsize}
\begin{align*}
&(69347771822290530959179953385244041759/201984628329,\\
& 1043901467465442412996262336398161817260855826874630788912/90777345556017483)\\
&(3460574163482908188427426712949356311/12038258961, \\
&-18110589319957094495275336848944697540739742502451835824/1320825734941959)\\
&(223827367575145652904326922554924756119/1504896849,\\
& 3348552250782159070948172554234175874539586500409115908848/58379463463257)\\
&(231186011876763383027626811894567567679/4259781289,\\
&3514369241948294947306662571663183913647384259899504483568/278023145389163),\\
&(24094394300976927534448792419082071/80982001, \\
&9710100599234190649985755744112823702319199233906064/728757026999),\\
&(4326551170759510514050992088776186591/12882477001, \\
&-17261150164745371973697619271361686469220145295462639664/1462174022090501),\\
&(27279101467693705413869829787600198711/54437755761, \\
&-45624789861369119683149199937634889173545547901847704624/12701362736400759),\\
&(22716106097239890740417034083, -3419512209640850893333977498177950568504084),\\
&(511456149321465121080387323, 2726580252438987666764744334964171072596),\\
&(948162104248407301294432275885974076879/1839096713689,\\
& 5772567839109568842182578374245288893337875076330512737232/2494059743625204637),\\
&(495179529275123438140025363, -4022507851813256992590035603118632514876),\\
&(1770287854471182905213974043, 61549688245299760953310275285220691708364).
\end{align*}
\end{scriptsize}
The regulator of these 12 points, as computed by SAGE, is $-2927257.43$. Since this is nonzero, we conclude that the above 12 points are independent.
 
Thus the generic rank of the family of curves given by \eqref{ecfam5} is at least 12.

\section{Concluding remarks}
While we have given several parametrized families of elliptic curves in Sections 2 and 3, the general method described in this paper can be used to obtain several other such families of curves.  In fact,  by solving certain symmetric diophantine systems, it may be possible to construct  parametrized families of elliptic curves with generic ranks $> 12$. Further, with efficient search techniques and greater computing power, the families of elliptic curves already obtained in this paper may yield several elliptic curves of fairly high rank. 

\begin{center}
\Large
Acknowledgments
\end{center}
 
I wish to  thank the Harish-Chandra Research Institute, Allahabad for providing me with all necessary facilities that have helped me to pursue my research work in mathematics. I also wish to thank Professor Andrew Bremner for doing some of the computations for this paper on the software MAGMA.

\noindent Postal Address: Ajai Choudhry, \\
\noindent \hspace{1.05 in} 13/4 A Clay Square, \\
\noindent \hspace{1.05 in} Lucknow-226001, INDIA.

\noindent E-mail address: ajaic203@yahoo.com


\begin{thebibliography}{99}


\bibitem{BMSW} Baur Bektemirov, Barry Mazur, William Stein, and Mark Watkins, {\it Average ranks of elliptic
curves: tension between data and conjecture}, Bull. Amer. Math. Soc. (N.S.) 44 (2007), 233--254.

\bibitem{ChW} A. Choudhry and  J. Wr\'oblewski, {\it Triads of integers with equal sums of squares and equal products and a related multigrade chain}, Acta Arithmetica,  178 (2017),  87--100.

\bibitem{Ca1} J. W. S. Cassels, {\it Diophantine equations with special reference to elliptic curves}, J. London Math. Soc. {\bf 41} (1966), 193--291.
\bibitem{Ca2} J. W. S. Cassels, Lectures on elliptic curves, Cambridge University Press, 1991. 

\bibitem{Di} L. E. Dickson, Introduction to the theory of numbers, Dover Publications, New York, 1957.

\bibitem{El} N. D. Elkies,  $\mathbb{Z}^{28}$ in  $E/\mathbb{Q}$, etc., 3 May 2006 e-mail available at NMBRTHRY@LISTSERV.NODAK.EDU

\bibitem{Fe} S. Fermigier, {\it Construction of high-rank elliptic curves over $\mathbb{Q}$ and $\mathbb{Q}(t)$ with non-trivial
2-torsion (extended abstract)}, Algorithmic number theory (Talence, 1996), 115--120, Lecture Notes in Comput. Sci.,
1122, Springer, Berlin, 1996.

\bibitem{Ki} S. Kihara, {\it On an elliptic curve over $\mathbb{Q}(t)$ of rank $\geq 14$}, Proc. Japan Acad. Ser. A, 77,  (2001), 50--51.

\bibitem{Me1} J.-F.  Mestre, {\it Courbes elliptiques de rang $≥11$ sur $Q(t)$}, C. R. Acad. Sci. Paris Sér. I Math. 313 (1991), 139--142.

\bibitem{Me2} J.-F.  Mestre, {\it Courbes elliptiques de rang $≥12$ sur $Q(t)$},  C. R. Acad. Sci. Paris Sér. I Math. 313 (1991), 171–-174.
 
\bibitem{Mo} L. J. Mordell, Diophantine equations, Academic Press, London, 1969. 

\bibitem{Na} K. Nagao, {\it An example of elliptic curve over $\mathbb{Q}(T)$ with rank $\geq 13$}, Proc. Japan Acad. Ser. A, 70,  (1994), 152--153.


\bibitem{PP} J. Park, B. Poonen, J, Voight and M. M. Wood, {\it A heuristic for boundedness of ranks of
elliptic curves}, available at arXiv:1602.01431 

\bibitem{Po} B. Poonen, {\it Heuristics for the arithmetic of elliptic curves}, available at  	arXiv:1711.10112 

\bibitem{SZ} S. Schmitt and H. G. Zimmer,  Elliptic Curves: A Computational Approach, Walter de Gruyter, Berlin 2003.

\bibitem{Sh} T. Shioda, {\it Construction of elliptic curves over $\mathbb{Q}(t)$ with high rank: a preview}, 
Proc. Japan Acad. Ser. A, 66 (1990), 57--60.

\bibitem{Si} J. H. Silverman, The arithmetic of elliptic curves, Second edition, Springer, Dordrecht, 2009.


\end{thebibliography}
	\end{document}